\documentclass[12pt,draftcls,onecolumn]{IEEEtran}
\usepackage{pdfsync}
\usepackage{graphicx}
\usepackage{amssymb}
\usepackage{amsfonts}
\usepackage{mathrsfs}
\usepackage{marvosym}
\usepackage[pdftex]{color}

\usepackage[dvips]{epsfig} \DeclareGraphicsExtensions{.ps,.eps}
\def\R{\mathbb{R}}

\def\det{{\rm det}\,}
\def\exp{{\rm exp}\,}

\def\tr{{\rm tr}\,}
\newcommand{\tp}{^{\top}}

\newcommand{\beq}{\begin{equation}}
\newcommand{\eeq}{\end{equation}}
\newcommand{\bea}{\begin{eqnarray}}
\newcommand{\eea}{\end{eqnarray}}

\newcommand{\bsea}{\begin{subeqnarray}}
\newcommand{\esea}{\end{subeqnarray}}
\newcommand{\nn}{\nonumber}

\newcommand{\proof}{\noindent {\it Proof. }}

\newcommand{\qed}{\hfill $\Box$ \vskip 2ex}

\def\bmat{\left[ \begin{array}}
\def\emat{\end{array} \right]}

\newcounter{acount}

\newtheorem{theorem}{Theorem}[section]

\newtheorem{lemma}[theorem]{Lemma}

\newtheorem{problem}[theorem]{Problem}

\begin{document}
\title{Matrix Completion by  the Principle of Parsimony}
\author{Augusto~Ferrante and Michele~Pavon\thanks{Work partially supported by the Italian Ministry for Education and Resarch (MIUR) under PRIN grant
``Identification and Robust Control of Industrial Systems", by the CPDA080209/08 and QFUTURE research grants of the University of Padova and by the Department of Information Engineering research project ``QUINTET".}\thanks{A. Ferrante is  with the
Dipartimento di Ingegneria dell'Informazione, Universit\`a di Padova,
via Gradenigo 6/B, 35131
Padova, Italy. {\tt\small  augusto@dei.unipd.it}}
\thanks{M. Pavon is with
the Dipartimento di Matematica Pura ed Applicata, Universit\`a di Padova, via
Trieste 63, 35131 Padova, Italy {\tt\small
pavon@math.unipd.it}}}

\markboth{DRAFT}{Shell \MakeLowercase{\textit{et al.}}: Bare Demo of IEEEtran.cls for Journals}
\maketitle

\begin{abstract}
Dempster's covariance selection method is extended first to general nonsingular matrices and then to full rank rectangular matrices. Dempster observed that his completion solved a maximum entropy problem. We show that our generalized completions are also solutions of a suitable entropy-like variational problem.
\end{abstract}


\IEEEpeerreviewmaketitle

\section{Dempster's covariance selection}
In the seminal paper \cite{Dempster-72}, Dempster introduced a general strategy for completing a partially specified covariance matrix.  Consider a zero-mean, multivariate Gaussian distribution with density
$$p(x)=(2\pi)^{-n/2}|\Sigma|^{-1/2}\exp\left\{-\frac{1}{2}x\tp\Sigma^{-1}x\right\},\quad x\in\R^n.
$$
Suppose that the elements $\{\sigma_{ij};1\le i\le j\le n, (i,j)\in \bar{{\cal I}}\}$ have been specified. How should $\Sigma$ be completed? Dempster resorts to a form of the {\em Principle of Parsimony} in parametric model fitting: As the elements $\sigma^{ij}$ of $\Sigma^{-1}$ appear as natural parameters of the model, one should set  $\sigma^{ij}$ to zero for $1\le i\le j\le n, (i,j)\not\in \bar{{\cal I}}$. Notice that $\sigma^{ij}=0$ has the probabilistic interpretation  that the $i$-th and $j$-th components of the Gaussian random vector  are {\em conditionally independent} given the other components.
This choice, which we name henceforth Dempster's Completion, may at first look less natural than setting the unspecified elements of $\Sigma$ to zero. It has nevertheless  considerable advantages compare to the latter, c.f. \cite[p.161]{Dempster-72}. In particular, Dempster  established the following far reaching result.
\begin{theorem}Assume that a symmetric, positive-definite completion of $\Sigma$ exists. Then there exists a unique Dempster's Completion $\Sigma^\circ$. This completion maximizes the (differential) entropy  
\begin{equation}\label{entropy}
H(p)=-\int_{\R^n}\log (p(x))p(x)dx=\frac{1}{2}\log(\det\Sigma)+\frac{1}{2}n\left(1+\log(2\pi)\right)
\end{equation}
among zero-mean Gaussian distributions having the prescribed elements $\{\sigma_{ij};1\le i\le j\le n, (i,j)\in \bar{{\cal I}}\}$.
\end{theorem} 
Thus, Dempster's Completion $\Sigma^\circ$ solves a {\em maximum entropy} problem, i.e. maximizes entropy under linear constraints.
Dempster's paper has generated a whole stream of research, see e.g \cite{BLW,SK,Dembo-M-S-89,HL} and references therein. In the meantime, matrix completion has become an important area of research with several new applications, where the completed matrix must have certain prescribed properties: For instance, it should be positive definite, it should be circulant, it should have a Toeplitz structure, it should have a prescribed low rank, etc. Motivation originates from problems in texture images modeling, recommender systems and networked sensors \cite{BGP,GUPTA,KLS1,CR,CP,KLS2,KMO,LV,CFPP}.

In this paper, we consider a totally unstructured version of Dempster's problem. A square matrix $\Sigma$ is partially specified and we seek to complete it according to the principle of parsimony. Besides the above mentioned applications, we are also motivated by the following problem. Suppose that we need to solve the linear system
\begin{equation}\Sigma X= B,
\end{equation}
where the matrix $B$ is given. Suppose that only the elements  $\{\sigma_{ij} \in \bar{{\cal I}}\}$ of $\Sigma$ could be estimated/determined, where $\bar{{\cal I}}$ is {\em any} subset of $ \{1,2,\dots n\}\times \{1,2,\dots n\}$. As a generic completion $\bar{\Sigma}$ is invertible, we can associate to such a completion the solution
\begin{equation}X_{\bar{\Sigma}}=\bar{\Sigma}^{-1} B.
\end{equation}  
We then identify as desirable, according to the principle of parsimony, completions $\bar{\Sigma}$ such that $\bar{\Sigma}^{-1}$ has a maximum number of zero entries. In this paper, we show that a family of such desirable completions are {\em generalized Dempster's completions} that can be characterized as critical points of a suitable variational problem. 

It may, at first, look hopeless to obtain an entropy-like variational characterization of Dempster's-like completions without positivity. Nevertheless, we {\em prove} in Lemmata \ref{derivata} and \ref{chainr} below that the contrained extremization of the determinant only involves the {\em positive part} of the matrix $\left(\Sigma\Sigma\tp\right)^{1/2}$. More precisely, only the {\em singular values} of $\Sigma$ come into play. Hence, such a variational characterization is possible and may be established even in the rectangular case.

The paper is outlined as follows. We discuss first the square case to facilitate the comparison with Dempster's classical results, see Sections \ref{prelim} and \ref{sec2} below. In Section \ref{secex}, we discuss two examples to illustrate the properties that our solutions may or may not enjoy. In Section \ref{rectangular}, we  generalize our results to rectangular, full rank matrices. The paper concludes with a discussion section comparing our approach to other matrix completion techniques and to other moment problems.

\section{The general completion problem}

Let ${\cal I}\subset \{1,2,\dots n\}\times \{1,2,\dots n\}$ and $\bar{{\cal I}}$ be the complementary subset.
To each $(i,j)\in{\cal I}$ we associate the unknown $x_{ij}$. Let $x$ be the vector, say $k$-dimensional, obtained by stacking the $x_{ij}$ one on top of the other.
We define as {\em partial matrix} a parametric family of matrices $\Sigma(x)$ whose entries
$[\Sigma(x)]_{i,j}=\sigma_{ij}$ are specified  for $(i,j)\in\bar{{\cal I}}$, while $[\Sigma(x)]_{i,j}=x_{ij}$ for $(i,j)\in{\cal I}$. Here, both $\sigma_{ij}$ and $x_{ij}$ take real values.
A {\em completion} of the partial matrix $\Sigma(x)$ is a matrix $\Sigma(\bar{x})$ where  $\bar{x}\in\R^k$. Notice that completions {\em always exist} as we are not requiring $\Sigma(\bar{x})$ to possess any further property.
If $\Sigma(x)$ is a partial matrix and $\cal I$ is the corresponding set of indices of the unspecified elements, we denote by ${\cal I}\tp$ the set of indices ${\cal I}\tp:=\{(j,i):\  (i,j)\in{\cal I}\}$.

Let  $\Sigma(x)$ be a square partial matrix of size $n$ and let $\cal I$ be
the corresponding set of indices of the unspecified entries.
Consider the following matrix completion problems:

\begin{problem}\label{problem1}
Find the nonsingular  completions $\Sigma(\bar{x})$ such that  
$[\Sigma(\bar{x})^{-1}]_{i,j}=0$ for all $(i,j)\in {\cal I}\tp$.
\end{problem}
\begin{remark}
Since in $\Sigma(x)$ we have $|{\cal I}|$ degrees of freedom (unknowns), we may generically expect
that the maximum number of entries that can be annihilated in $\Sigma(x)^{-1}$ is precisely $|{\cal I}|$. Indeed, $m$ zeros in $\Sigma(x)^{-1}$ is equivalent to finding a solution to a system of $m$ polynomial equations in $|{\cal I}|$ unknowns on the complement of the set where the determinant vanishes (in Section \ref{secex}, we provide a non generic example where there is a completion with more than $|{\cal I}|$ zeros in the inverse). 
\end{remark}
\begin{problem}\label{problem2}
Find the nonsingular completions $\Sigma(\bar{x})$ that extremize  
$\det[\Sigma(x)]$.
\end{problem}
\begin{remark} Notice that when a covariance matrix $\Sigma=\Sigma\tp>0$ is sought, Problem \ref{problem2} reduces to the maximum entropy problem solved by Dempster's Completion. This follows from the fact that the entropy, in the Gaussian case, differs from $(1/2)\log\det\Sigma$ by a constant, the monotonicity of the logarithm and strict concavity of the entropy. Thus, Problem \ref{problem2} appears as a legitimate generalization of Dempster's classical completion method. 
\end{remark}

The main result of  this paper consists in showing  that Problems \ref{problem1} and \ref{problem2} have the same set of solutions.
\begin{theorem}\label{mainres}
$\Sigma(\bar{x})$ solves {\em Problem \ref{problem1}} if and only if it solves {\em Problem \ref{problem2}}.
\end{theorem}
\begin{remark}
It is apparent that solutions to Problems \ref{problem1} and \ref{problem2} may not exist, but when they do there may be many. For instance, consider
$$\Sigma(x)=\left(
\begin{array}{cc}
 1 & x_{12} \\
 3 & x_{22}
\end{array}
\right).
$$
Then, Problems \ref{problem1} and \ref{problem2} are not solvable. 
Actually, whenever all the unknowns $x_{ij}$ are in the same row or in the same column, $\det(\Sigma(x))$ is linear in $x_{ij}$ and hence it does not have critical points.
Two examples where there exist multiple solutions are provided in Section \ref{secex}.
\end{remark}

\section{Some preliminary results}\label{prelim}
We collect below some lemmata that are needed to prove Theorem \ref{mainres}.
\begin{lemma}\label{derivata}
Problem \ref{problem2} is equivalent to the following: Find a nonsingular completion $\Sigma(\bar{x})$ that extremizes 
$J(\Sigma(x)):=\log|\det[\Sigma(x)]|$.
\end{lemma}
\proof
Compute the gradient
\beq
\frac{\partial}{\partial x} [\log|\det[\Sigma(x)]|]=\frac{1}{|\det[\Sigma(x)]|}\frac{|\det[\Sigma(x)]|}{\det[\Sigma(x)]}\frac{\partial}{\partial x}[\det[\Sigma(x)]]=\frac{1}{\det[\Sigma(x)]}\frac{\partial}{\partial x}[\det[\Sigma(x)]]
\eeq
Now the statement follows by observing that we are restricting attention to nonsingular completions.
\qed

Denote by $D[J(\Sigma);\delta\Sigma]$ the directional derivative of $J$ in direction $\delta\Sigma\in\R^{n\times n}$:
\beq
D[J(\Sigma);\delta\Sigma]=\lim_{\varepsilon\rightarrow 0}\frac{\log|\det[\Sigma+\varepsilon\delta\Sigma]|-\log|\det[\Sigma]|}{\varepsilon}.
\eeq
We have the following result.
\begin{lemma}\label{chainr}
Let $J(\Sigma(x))=\log|\det[\Sigma(x)]|$ as in the previous lemma. If $\Sigma$ is nonsingular then, for any $\delta\Sigma\in\R^{n\times n}$,
\beq
D[J(\Sigma);\delta\Sigma]=\tr[\Sigma^{-1}\delta\Sigma].
\eeq
\end{lemma}
\proof
Let
\beq \label{polardec}
P(\Sigma):=\left(\Sigma\Sigma\tp\right)^{1/2},\quad U(\Sigma):=P(\Sigma)^{-1}\Sigma.
\eeq
Observe that $\Sigma =P(\Sigma)U(\Sigma)$ is the polar decomposition of $\Sigma$. Similarly, we have
$
\Sigma_\varepsilon:=\Sigma+\varepsilon\delta\Sigma=P(\Sigma_\varepsilon)U(\Sigma_\varepsilon).
$
Consider now the Taylor expansion
\beq
P(\Sigma_\varepsilon)=P(\Sigma) + \varepsilon\delta P + o(\varepsilon^2),
\eeq
where $\delta P:= D[P(\Sigma);\delta\Sigma]$. In view of (\ref {polardec}), the latter directional derivative may be expressed as
\beq \label{deltaP}
\delta P= D[P(\Sigma);\delta\Sigma]=D[\Sigma U\tp(\Sigma);\delta\Sigma]=\delta\Sigma U\tp(\Sigma)+\Sigma D[U\tp(\Sigma);\delta\Sigma].
\eeq

Notice now that 
\beq\label{fatto1}
\log|\det[\Sigma]|=\log[\det[P(\Sigma)]].
\eeq
Moreover, if $Q=Q\tp>0$, the following expression holds \cite{BLW,GEORGIOU_RELATIVEENTROPY}:
\beq\label{derdp}
D[\log[\det(Q)];\delta Q]=\tr[Q^{-1}\delta Q].
\eeq
From (\ref{fatto1}), (\ref{derdp}) and  (\ref{deltaP}), we now get:
\bea\nn
D[J(\Sigma);\delta\Sigma]&=&\tr[P(\Sigma)^{-1}\delta P]\\
\nn
&=&\tr[P(\Sigma)^{-1}\delta\Sigma U\tp(\Sigma)]+\tr[P(\Sigma)^{-1}\Sigma D[U\tp(\Sigma);\delta\Sigma]]\\
&=& \tr[\Sigma^{-1}\delta\Sigma] +\tr[U(\Sigma)  D[U\tp(\Sigma);\delta\Sigma]].
\eea
The result now follows by observing that $\tr[U(\Sigma)  D[U\tp(\Sigma);\delta\Sigma]]=0$.
Indeed,
\bea\nn
0&=&  \tr[D[U(\Sigma)U\tp(\Sigma);\delta\Sigma]]\\
\nn
&=&\tr[U(\Sigma)  D[U\tp(\Sigma);\delta\Sigma]]+\tr[D[U(\Sigma);\delta\Sigma]U\tp(\Sigma)]\\
&=& 2\;\tr[U(\Sigma)  D[U\tp(\Sigma);\delta\Sigma]].\label{tr=0}
\eea
\qed

\begin{lemma}\label{sporth}
Consider the space $\R^{n\times n}$ endowed with the inner product $\langle M_1,M_2\rangle:=\tr[M_1\tp M_2]$.
Let ${\cal M}$ be the subspace of $\R^{n\times n}$ consisting of the matrices whose entries in position $(i,j)\in \bar{{\cal I}}$ are zero.
Let $M\in{\cal M}^\perp$. Then, $[M]_{i,j}=0$ for all $(i,j)\in{\cal I}$.
\end{lemma}
\proof
Denote by $e_i$ the $i$-th canonical vector in $\R^n$. Clearly, for any $(i,j)\in{\cal I}$, $e_ie_j\tp\in{\cal M}$.
Thus, if $M\in{\cal M}^\perp$, for any $(i,j)\in{\cal I}$, we have $0=\tr[(e_ie_j\tp)\tp M]=
\tr[e_je_i\tp M]=e_i\tp Me_j=[M]_{ij}$.
\qed

\section{Proof of Theorem \ref{mainres}}\label{sec2}

We are now ready to prove our main result.

{\em Proof of Theorem \ref{mainres}.}
In view of Lemma \ref{derivata}, Problem \ref{problem2} is equivalent to 
the following variational problem:
 \beq
 {\rm extremize}\left\{J(\Sigma):\ e\tp_i\Sigma e_j=\sigma_{ij},\ (i,j)\in\bar{{\cal I}}\right\},
 \eeq
where, as before, $J(\Sigma(x))=\log|\det[\Sigma(x)]|$. The corresponding Lagrangian is
\beq\label{Lagrangian}
{\cal L}(\Sigma)= J(\Sigma)+\sum_{(i,j)\in\bar{{\cal I}}} \lambda_{ij}( e\tp_i\Sigma e_j-\sigma_{ij})
 \eeq
whose {\em unconstrained} extremization is obtained by annihilating the directional derivative of ${\cal L}$ in any direction $\delta\Sigma\in\R^{n\times n}$. In view of Lemma \ref{chainr},
this yields:
\beq
\tr\left[\left(\Sigma^{-1}-\sum_{(i,j)\in\bar{{\cal I}}}\lambda_{ij}e_je\tp_i\right) \delta\Sigma\right]=0,\quad \forall \ \delta\Sigma\in\R^{n\times n},
\eeq
or, equivalently,
\beq\label{formofoptinv}
\Sigma^{-1}=\sum_{(i,j)\in\bar{{\cal I}}}\lambda_{ij}e_je\tp_i.
\eeq
It follows, in particular, that the inverse of any nonsingular critical point $\Sigma$ of the Lagrangian (\ref{Lagrangian}) has zeros in positions $(i,j)\in{\cal I}\tp$. Moreover, if we can find  $\lambda^\circ_{ij}$  such that the matrix $\displaystyle{\sum_{(i,j)\in\bar{{\cal I}}}\lambda^\circ_{ij}e_je\tp_i}$ is nonsingular and
\beq\label{formsopt}\Sigma^\circ:={\displaystyle{\left(\sum_{(i,j)\in\bar{{\cal I}}}\lambda^\circ_{ij}e_je\tp_i\right)^{-1}}}
\eeq
 satisfies
\beq\label{condopt}
e_i\tp\Sigma^\circ e_j=\sigma_{ij},\quad (i,j)\in\bar{{\cal I}},
\eeq
then $\Sigma^\circ $ is indeed a solution of Problem \ref{problem2}. 

Assume now that $\Sigma(\bar{x})$ solves Problem \ref{problem1}. Then,
$[\Sigma(\bar{x})]^{-1}$  has
the form  (\ref{formofoptinv}). Moreover, $\Sigma(\bar{x})$ satisfies (\ref{condopt}). Hence, it solves 
Problem \ref{problem2}.
 
Conversely,  let $\Sigma(\bar{x})$ solve Problem \ref{problem2}.
This is equivalent to 
\beq
 D[J(\Sigma),\delta\Sigma]_{{\big |}\Sigma=\Sigma(\bar{x})}=\tr[\Sigma(\bar{x})^{-1}\delta\Sigma]=\langle \Sigma(\bar{x})^{-\top},\delta\Sigma\rangle=0,\qquad
 \forall\ \delta\Sigma\in{\cal M},
 \eeq
 where ${\cal M}$ is the subspace of $\R^{n\times n}$
 (defined in Lemma \ref{sporth}) consisting of  matrices whose entries in position $(i,j)\in \bar{{\cal I}}$ are zero.
 Thus, $\Sigma(\bar{x})^{-\top}\in{\cal M}^\perp$. By Lemma \ref{sporth}, $[\Sigma(\bar{x})^{-1}]_{i,j}=0$
 for all $(i,j)\in {\cal I}\tp$, namely $\Sigma(\bar{x})$ solves Problem \ref{problem1}.
 \qed

\section{Two illustrative examples}\label{secex}
{\bf Example 1.} The following example shows that in some pathological situations it is indeed possible to complete a partial matrix in such a way that the completion is symmetric and positive definite and has a larger number of vanishing entries than the Dempster completion.
Consider the matrix

$$
\Sigma(x)=\left(
\begin{array}{cccc}
 \frac{120}{929} & \frac{4}{929} & -\frac{15}{929} & x \\
 \frac{4}{929} & \frac{124}{929} & -\frac{1}{1858} & -\frac{63}{1858} \\
 -\frac{15}{929} & -\frac{1}{1858} & \frac{118}{929} & \frac{2}{929} \\
 x & -\frac{63}{1858} & \frac{2}{929} & \frac{126}{929}
\end{array}
\right)
$$

Dempster's completion corresponds to $x=x_d=-(79/58527)$. The associated matrix
$\Sigma(x_d)$ has the following inverse:

$$\Sigma(x_d)^{-1}=\left(
\begin{array}{cccc}
 \frac{63}{8} & -\frac{1}{4} & 1 & 0 \\
 -\frac{1}{4} & \frac{1009}{126} & -\frac{2}{63} & 2 \\
 1 & -\frac{2}{63} & \frac{4033}{504} & -\frac{1}{8} \\
 0 & 2 & -\frac{1}{8} & \frac{63}{8}
\end{array}
\right)$$

If, on the other hand, we pick  $x=x_m=-(16/929)$, the associated matrix
$\Sigma(x_m)$ has the following inverse:
$$
\Sigma(x_m)^{-1}=\left(
\begin{array}{cccc}
 8 & 0 & 1 & 1 \\
 0 & 8 & 0 & 2 \\
 1 & 0 & 8 & 0 \\
 1 & 2 & 0 & 8
\end{array}
\right)$$
which has six vanishing entries.
Notice that, as for the Dempster's Completion,  this extension is  symmetric and positive definite.

\noindent
{\bf Example 2.} Let 
$$\Sigma(x)=\left(
\begin{array}{cccc}
 5 & 2 & x & 1 \\
 2 & 5 & 2 & y \\
 w & 2 & 5 & 2 \\
 1 & z & 2 & 5
\end{array}
\right).$$
Notice that the specified entries of $\Sigma(x)$ are compatible with symmetry and  the Toeplitz structure.
By extremizing $\det[\Sigma(x)]$, we obtain seven real matrices completing $\Sigma(x)$ to a nonsingular matrix: 
$$\Sigma(\bar{x}_1)=\left(
\begin{array}{cccc}
 5 & 2 & -6 & 1 \\
 2 & 5 & 2 & -6 \\
 -6 & 2 & 5 & 2 \\
 1 & -6 & 2 & 5
\end{array}
\right),\quad
 \Sigma(\bar{x}_2)=\left(
\begin{array}{cccc}
 5 & 2 & -\frac{19}{5} & 1 \\
 2 & 5 & 2 & 5 \\
 -\frac{19}{5} & 2 & 5 & 2 \\
 1 & 5 & 2 & 5
\end{array}
\right),$$
$$
\Sigma(\bar{x}_3)= \left(
\begin{array}{cccc}
 5 & 2 & 1 & 1 \\
 2 & 5 & 2 & 1 \\
 1 & 2 & 5 & 2 \\
 1 & 1 & 2 & 5
\end{array}
\right),\
 \Sigma(\bar{x}_4)=\left(
\begin{array}{cccc}
 5 & 2 & 5 & 1 \\
 2 & 5 & 2 & -\frac{19}{5} \\
 5 & 2 & 5 & 2 \\
 1 & -\frac{19}{5} & 2 & 5
\end{array}
\right),\
\Sigma(\bar{x}_5)= \left(
\begin{array}{cccc}
 5 & 2 & 5 & 1 \\
 2 & 5 & 2 & 5 \\
 5 & 2 & 5 & 2 \\
 1 & 5 & 2 & 5
\end{array}
\right),
$$
$$\Sigma(\bar{x}_6)= \left(
\begin{array}{cccc}
 5 & 2 & -\frac{3}{2} \left(\sqrt{13}-5\right) & 1 \\
 2 & 5 & 2 & -\frac{3}{2} \left(\sqrt{13}-5\right) \\
 \frac{3}{2} \left(5+\sqrt{13}\right) & 2 & 5 & 2 \\
 1 & \frac{3}{2} \left(5+\sqrt{13}\right) & 2 & 5
\end{array}
\right),$$
$$
\Sigma(\bar{x}_7)= \left(
\begin{array}{cccc}
 5 & 2 & \frac{3}{2} \left(5+\sqrt{13}\right) & 1 \\
 2 & 5 & 2 & \frac{3}{2} \left(5+\sqrt{13}\right) \\
 -\frac{3}{2} \left(\sqrt{13}-5\right) & 2 & 5 & 2 \\
 1 & -\frac{3}{2} \left(\sqrt{13}-5\right) & 2 & 5
\end{array}
\right).
$$

All of these completions have inverse with zeros in positions $(1,3), (2,4), (3,1)$ and $(4,2)$.
Notice that only $\Sigma(\bar{x}_1)$, $\Sigma(\bar{x}_3)$, and $\Sigma(\bar{x}_5)$
are symmetric and have a Toeplitz structure. Among these, only $\Sigma(\bar{x}_3)$ is also positive definite
($\Sigma(\bar{x}_3)$ is indeed the Dempster Completion).
$\Sigma(\bar{x}_2)$, $\Sigma(\bar{x}_4)$ are symmetric but do not have a Toeplitz structure and $\Sigma(\bar{x}_6)$, $\Sigma(\bar{x}_7)$ have a Toeplitz structure but they are not symmetric.

We observe that all of the $5$ symmetric completions are also solutions of the problem of extremizing 
$\det\left(
\begin{array}{cccc}
 5 & 2 & x & 1 \\
 2 & 5 & 2 & y \\
 x & 2 & 5 & 2 \\
 1 & y & 2 & 5
\end{array}
\right)$.

Similarly, all of the $5$ Toeplitz completions are also solutions of the problem of extremizing 
$\det\left(
\begin{array}{cccc}
 5 & 2 & x & 1 \\
 2 & 5 & 2 & x \\
 y & 2 & 5 & 2 \\
 1 & y & 2 & 5
\end{array}
\right)$.

This example shows that even if the constraints are compatible with symmetry or other matrix properties, there may exist extremizing completions that do not preserve these features.

\section{The case of rectangular matrices}\label{rectangular}

Next, we extend the results obtained in the previous sections to the general case of possibly non-square matrices $\Sigma(x)\in\R^{n\times p}$.
We assume that $p\geq n$ (the case  $p\leq n$ can be dealt with in a dual fashion). 
As before, let ${\cal I}\subset \{1,2,\dots n\}\times \{1,2,\dots p\}$, and $\bar{{\cal I}}$ be the complementary subset.
To each $(i,j)\in{\cal I}$ we associate the unknown $x_{ij}$. Let $x$ be the $k$-dimensional vector obtained by stacking the $x_{ij}$ one on top of the other.
Define as before a {\em partial matrix} to be a parametric family of matrices $\Sigma(x)$ whose entries
$[\Sigma(x)]_{i,j}=\sigma_{ij}$ are specified  for $(i,j)\in\bar{{\cal I}}$, while $[\Sigma(x)]_{i,j}=x_{ij}$ for $(i,j)\in{\cal I}$. Again, $\sigma_{ij}$ and $x_{ij}$ take real values.
Consider the following matrix completion problems:

\begin{problem}\label{problem1-gen}
Find full row rank  completions $\Sigma(\bar{x})$ such that the corresponding  Moore-Penrose pseudo-inverse  $\Sigma(\bar{x})^{\sharp}$ satisfies
$[\Sigma(\bar{x})^{\sharp}]_{i,j}=0$ for all $(i,j)\in {\cal I}\tp$.
\end{problem}
Notice that, since $\Sigma(\bar{x})$ is full row rank, the Moore-Penrose pseudo-inverse is also a right-inverse and is explicity given by
\begin{equation}\label{MPPI}
\Sigma(\bar{x})^{\sharp}=\Sigma(\bar{x})\tp\left(\Sigma(\bar{x})\Sigma(\bar{x})\tp\right)^{-1}.
\end{equation}

\begin{problem}\label{problem2-gen}
Find full row rank   completions $\Sigma(\bar{x})$ that extremize  
$\det\left[\left(\Sigma(x)\Sigma(x)\tp\right)^{1/2}\right]$.
\end{problem}

\begin{remark}
The form of the index in Problem \ref{problem2-gen} is inspired by the fact, established in Lemmata \ref{derivata} and \ref{chainr}, that in the square case the variational analysis only depends on the positive part $P(\Sigma(x))=[\Sigma(x)\Sigma(x)\tp]^{1/2}$ of $\Sigma(x)$. Actually, only the {\em singular values} of $\Sigma(x)$ come into play.
\end{remark}

The main result of  this section consists in showing  that Problems \ref{problem1-gen} and \ref{problem2-gen} are equivalent.
\begin{theorem}\label{mainres-gen}
$\Sigma(\bar{x})$ solves {\em Problem \ref{problem1-gen}} if and only if it solves {\em Problem \ref{problem2-gen}}.
\end{theorem}

We first notice that 
Problem \ref{problem2-gen} is equivalent to the following: Find a full row rank   completion $\Sigma(\bar{x})$ that extremizes 
\beq\label{JNS}
J(\Sigma(x)):=\log\left[\det\left[\left(\Sigma(x)\Sigma(x)\tp\right)^{1/2}\right]\right].
\eeq

Denote by $D[J(\Sigma);\delta\Sigma]$ the directional derivative of $J$ in direction $\delta\Sigma\in\R^{n\times p}$:
\beq
D[J(\Sigma);\delta\Sigma]=\lim_{\varepsilon\rightarrow 0}\frac{J(\Sigma+\varepsilon\delta\Sigma)-J(\Sigma)}{\varepsilon}.
\eeq
We have
\begin{lemma}\label{chainr-gen}
If $\Sigma$ is full row rank   then, for any $\delta\Sigma\in\R^{n\times p}$,
\beq
D[J(\Sigma);\delta\Sigma]=\tr[\Sigma^{\sharp}\delta\Sigma].
\eeq
\end{lemma}
\proof
Let
\beq \label{polardec-gen}
P(\Sigma):=(\Sigma\Sigma\tp)^{1/2},\quad U(\Sigma):=P(\Sigma)^{-1}\Sigma.
\eeq
Observe that $\Sigma =P(\Sigma)U(\Sigma)$ is the generalized polar decomposition of $\Sigma$. In particular, $U(\Sigma)$ is a matrix whose rows are orthonormal: $U(\Sigma)U(\Sigma)\tp=I$. From (\ref{MPPI}), we get the following representation for the Moore-Penrose pseudo-inverse of $\Sigma$
\begin{equation}\label{RGI}
\Sigma^{\sharp}=U(\Sigma)\tp P(\Sigma)^{-1}.
\end{equation}
Similarly, we have
$
\Sigma_\varepsilon:=\Sigma+\varepsilon\delta\Sigma=P(\Sigma_\varepsilon)U(\Sigma_\varepsilon).
$
Consider now the Taylor expansion
\beq
P(\Sigma_\varepsilon)=P(\Sigma) + \varepsilon\delta P + o(\varepsilon^2),
\eeq
where $\delta P:= D[P(\Sigma);\delta\Sigma]$. In view of (\ref {polardec-gen}), the latter directional derivative may be expressed as
\beq \label{deltaP-gen}
\delta P= D[P(\Sigma);\delta\Sigma]=D[\Sigma U\tp(\Sigma);\delta\Sigma]=\delta\Sigma U\tp(\Sigma)+\Sigma D[U\tp(\Sigma);\delta\Sigma].
\eeq

From  (\ref{derdp}), (\ref{deltaP-gen}) and (\ref{RGI}), using the cyclic property of the trace, we now get:
\bea\nn
D[J(\Sigma);\delta\Sigma]&=&\tr[P(\Sigma)^{-1}\delta P]\\
\nn
&=&\tr[P(\Sigma)^{-1}\delta\Sigma U\tp(\Sigma)]+\tr[P(\Sigma)^{-1}\Sigma D[U\tp(\Sigma);\delta\Sigma]]\\
&=& \tr[\Sigma^{\sharp}\delta\Sigma] +\tr[U(\Sigma)  D[U\tp(\Sigma);\delta\Sigma]].
\eea
The result now follows by observing that $\tr[U(\Sigma)  D[U\tp(\Sigma);\delta\Sigma]]=0$ as in (\ref{tr=0}).
\qed

The following result is a simple generalization of Lemma \ref{sporth}. 
\begin{lemma}\label{sporth-gen}
Consider the space $\R^{n\times p}$ endowed with the inner product $\langle M_1,M_2\rangle:=\tr[M_1\tp M_2]$.
Let ${\cal M}$ be the subspace of $\R^{n\times p}$ consisting of the matrices whose entries in position $(i,j)\in \bar{{\cal I}}$ are zero.
Let $M\in{\cal M}^\perp$. Then, $[M]_{i,j}=0$ for all $(i,j)\in{\cal I}$.
\end{lemma}

{\em Proof of Theorem \ref{mainres-gen}.}
Problem \ref{problem2-gen} is equivalent to 
the following variational problem
 \beq
 {\rm extremize}\left\{J(\Sigma):\ e\tp_i\Sigma e_j=\sigma_{ij},\ (i,j)\in\bar{{\cal I}}\right\},
 \eeq
 where $J$ is given by (\ref{JNS}).The corresponding Lagrangian is
\beq\label{Lagrangian-gen}
{\cal L}(\Sigma)= J(\Sigma)+\sum_{(i,j)\in\bar{{\cal I}}} \lambda_{ij}( e\tp_i\Sigma e_j-\sigma_{ij})
 \eeq
whose {\em unconstrained} extremization is obtained by annihilating the directional derivative of ${\cal L}$ in any direction $\delta\Sigma$. In view of Lemma \ref{chainr-gen},
we get :
\beq
\tr\left[\left(\Sigma^{\sharp}-\sum_{(i,j)\in\bar{{\cal I}}}\lambda_{ij}e_je\tp_i\right) \delta\Sigma\right]=0,\quad \forall \ \delta\Sigma\in\R^{n\times p},
\eeq
or, equivalently,
\beq\label{formofoptinv-gen}
\Sigma^{\sharp}=\sum_{(i,j)\in\bar{{\cal I}}}\lambda_{ij}e_je\tp_i.
\eeq
It follows, in particular, that the Moore-Penrose pseudo-inverse of any full row-rank critical point $\Sigma$ of the Lagrangian (\ref{Lagrangian}) has zeros in positions $(i,j)\in{\cal I}\tp$. Moreover, if we can find  $\lambda^\circ_{ij}$  such that the matrix $\displaystyle{\sum_{(i,j)\in\bar{{\cal I}}}\lambda^\circ_{ij}e_je\tp_i}$ is full column-rank and
\beq\label{formsopt-gen}\Sigma^\circ:={\displaystyle{\left(\sum_{(i,j)\in\bar{{\cal I}}}\lambda^\circ_{ij}e_je\tp_i\right)^{\sharp}}}
\eeq
 satisfies
\beq\label{condopt-gen}
e_i\tp\Sigma^\circ e_j=\sigma_{ij},\quad (i,j)\in\bar{{\cal I}},
\eeq
then $\Sigma^\circ $ is indeed a solution of Problem \ref{problem2-gen}. 

Assume now that $\Sigma(\bar{x})$ solves Problem \ref{problem1-gen}. Then, $[\Sigma(\bar{x})]$ is full row-rank and
$[\Sigma(\bar{x})]^{\sharp}$  has
the form  (\ref{formofoptinv-gen}). Moreover, $\Sigma(\bar{x})$ satisfies (\ref{condopt-gen}). Hence, it solves 
Problem \ref{problem2-gen}.
 
Conversely,  let $\Sigma(\bar{x})$ solve Problem \ref{problem2-gen}.
This is equivalent to 
\beq
 D[J(\Sigma),\delta\Sigma]_{{\big |}\Sigma=\Sigma(\bar{x})}=\tr[\Sigma(\bar{x})^{\sharp}\delta\Sigma]=\langle \left[\Sigma(\bar{x})^{\sharp}\right]^{\top},\delta\Sigma\rangle=0,\qquad
 \forall\ \delta\Sigma\in{\cal M},
 \eeq
 where ${\cal M}$ is the subspace of $\R^{n\times p}$
 (defined in Lemma \ref{sporth-gen}) consisting of   matrices whose entries in position $(i,j)\in \bar{{\cal I}}$ are zero.
 Thus, $\left[\Sigma(\bar{x})^{\sharp}\right]^{\top}\in{\cal M}^\perp$. By Lemma \ref{sporth-gen}, $[\Sigma(\bar{x})^{\sharp}]_{i,j}=0$
 for all $(i,j)\in {\cal I}\tp$, namely $\Sigma(\bar{x})$ solves Problem \ref{problem1-gen}.
 \qed

We outline again the significance of the above result for the solution of systems of linear equations. Consider
\begin{equation}X\Sigma = B,
\end{equation}
where $B$ is given. Suppose that only certain elements of $\Sigma$ could be estimated/determined. If the completion $\bar{\Sigma}$ has full row rank, we can associate to it the solution
\begin{equation}X_{\bar{\Sigma}}=B\bar{\Sigma}^{\sharp}.
\end{equation}  
Again we identify as desirable, according to the principle of parsimony, completions $\bar{\Sigma}$ such that $\bar{\Sigma}^{\sharp}$ has a maximum number of zero entries. 

\section{Discussion}
The nuclear norm ({\em sum of singular values}) of a matrix is often used in convex
heuristics for rank minimization problems in control, signal processing, and statistics. It has been employed in a series of recent papers on matrix completion, see \cite{CR,CP,KMO,LV} and references therein. The renewed interest in this metric has both theoretical and practical reasons as argued in the above mentioned papers. Variational problems involving the {\em sum of the logarithm of the singular values} (the logarithm of the determinant in the covariance case), such as those presented in this paper, occupy a somewhat complementary place. Indeed, as we have shown above, they lead to constraints on the (pseudo-)inverse of $\Sigma$.   Moreover, in the case when $\Sigma$ is a covariance matrix, Dempster's Completion $\Sigma^\circ$ {\em maximizes} entropy, namely the sum of the logarithm of the singular values (eigenvalues), whereas in  \cite{CR,CP,KMO,LV} the sum of the singular values is {\em minimized}.

Our variational problems appear close in spirit to Janes \cite{Jaynes57,Jaynes82}, where, following in the footsteps of Boltzmann (1877), Schr\"{o}dinger \cite{S}, Cram\'er  \cite{C}, Sanov  \cite{SANOV}, etc., and followed by such coryphaei as Dempster himself \cite{Dempster-72}, Akaike \cite{A1}, Burg \cite{B,BLW}, etc., he promoted maximum entropy methods to general inference methods. It might actually be worthwhile to quote an illuminating passage from the introduction of \cite{Jaynes82} which deals with spectral analysis: ``There are many different spectral analysis problems, corresponding to different kinds of prior information about the phenomenon being observed, different kinds of data, different kinds of perturbing noise, and different objectives. It is, therefore, quite meaningless to pass judgment on the merits of any proposed method unless one specifies clearly: ``In what class of problems is this method intended to be used?" Most of the current confusion on these questions is, in the writer's opinion, the direct result of failure to define the problem explicitly enough." We feel that these considerations apply equally well to the matrix completion problem. 

In this paper, we have chosen to discuss a very general completion problem where no further requirement is imposed on the solution matrix. As soon as the solution is required to feature some properties, such as being  positive definite, (with the possible additional constraints of being Toeplitz, circulant,etc.) the {\em existence} of matrices having the prescribed elements and properties becomes an issue. When existence is guaranteed, it should be apparent that our variational analysis can be readily  adapted to these more structured problems. 

We finally want to observe that completing a matrix so that it enjoys certain properties \cite{Dempster-72,BLW,SK,Dembo-M-S-89,BGP,GUPTA,KLS1,HL,CR,CP,KLS2,KMO,LV,CFPP} may be viewed as a {\em generalized moment problem}. These are problems where a function (a measure, a matrix, etc.) is sought satisfying certain given moment constraints as in the classical moment problem \cite{A,KN}, but also enjoying further properties: These may take several different forms. We mention the important case of bounds on the complexity, such as a bound on the degree of the rational solution, for applications in communications and control engineering, cf. e.g. \cite{GEORGIOU_REALIZATION,BGuL,BYRNES_GEORGIOU_LINDQUIST_THREE,BYRNES_GUSEV_LINDQUIST_FROMFINITE,GEORGIOU_LINDQUIST_KULLBACKLEIBLER,GEORGIOU_RELATIVEENTROPY,FERRANTE_PAVON_RAMPONI_HELLINGERVSKULLBACK,RAMPONI_FERRANTE_PAVON_GLOBALLYCONVERGENT}.

\end{document}